\documentclass{article}

\usepackage{amsmath,amsfonts,amsthm,amssymb}
\usepackage{url}

\theoremstyle{plain}
\newtheorem{theorem}{Theorem}
\newtheorem{corollary}[theorem]{Corollary}
\newtheorem{proposition}[theorem]{Proposition}

\theoremstyle{definition}
\newtheorem{remark}[theorem]{Remark}
\newtheorem{example}[theorem]{Example}
\newtheorem{definition}[theorem]{Definition}

\begin{document}

\title{The Natural Logarithm on Time Scales}

\author{Dorota Mozyrska\thanks{Supported by
        Bia\l ystok Technical University grant S/WI/1/08.}\\
        \url{admoz@w.tkb.pl}\\[0.3cm]
        Faculty of Computer Science\\
        Bia{\l}ystok Technical University\\
        15-351 Bia\l ystok, Poland
    \and
        Delfim F. M. Torres\thanks{Supported by
        the R\&D unit CEOC, via FCT
        and the EC fund FEDER/POCI 2010.}\\
        \url{delfim@ua.pt}\\[0.3cm]
        Department of Mathematics\\
        University of Aveiro\\
        3810-193 Aveiro, Portugal}

\date{September 2008}

\maketitle


\begin{abstract}
We define an  appropriate logarithm
function on time scales and present its main properties.
This gives answer to a question posed by M.~Bohner in
[J. Difference Equ. Appl. {\bf 11} (2005), no.~15, 1305--1306].

\bigskip


\medskip

\end{abstract}

\maketitle


\section{Introduction}

The following open problem is posed in \cite{OpenQuest}:
\begin{quotation}
\emph{Define a ``nice'' logarithm function on time scales.}
\end{quotation}
By ``nice'' it is understood a function $L(\cdot)$ characterized
by the equation $L(x y) = L(x) + L(y)$ and such that
$L(t) = \ln(t)$ when $\mathbb{T} = \mathbb{R}$.
In this short note we give a simple answer to this problem.


\section{The definition}

An introduction to time scales can be found in \cite{book:ts}.
Throughout the text we assume $\mathbb{T}$ to be a time scale with at least two positive points, one of them being always one:
$1 \in \mathbb{T}$, there exists at least one $t \in \mathbb{T}$ such that $0 < t \ne 1$. We define the \emph{natural logarithm} function on the time scale $\mathbb{T}$ by
\begin{equation}
\label{eq:def}
L_{\mathbb{T}}(t) :=
\int_1^t \frac{1}{\tau} \Delta\tau \, , \quad
t \in \mathbb{T} \cap (0,+\infty) \, .
\end{equation}
The motivation for our definition \eqref{eq:def} is clear:
$L_{\mathbb{T}}(\cdot)$ is a delta differentiable function
satisfying
\begin{equation}
\label{eq:trv:propt}
L_{\mathbb{T}}^{\Delta}(t) = \frac{1}{t} \, , \quad
L_{\mathbb{T}}(1) = 0 \, , \quad \text{and }
L_{\mathbb{R}}(t) = \ln(t) \, .
\end{equation}
Since $L_{\mathbb{T}}^{\Delta}(t) > 0$ for every
$t \in \mathbb{T}^{\kappa} \cap (0,+\infty)$,
$L_{\mathbb{T}}(\cdot)$ is an increasing
and continuous function.
Moreover, for $t \in \mathbb{T} \cap (0,+\infty)$,
$L_{\mathbb{T}}(t)\leq t-1$;
if $t < 1$, then $L_{\mathbb{T}}(t) < 0$;
if $t > 1$, then $L_{\mathbb{T}}(t) > 0$.
Directly from the definition we have also that $L_{\mathbb{T}}(\sigma(t))=L_{\mathbb{T}}(t)
+\mu(t)L_{\mathbb{T}}^{\Delta}(t)=L_{\mathbb{T}}(t)+\frac{\mu(t)}{t}$, where $\mu(t)=\sigma(t)-t$.

We remark that definition \eqref{eq:def} is different from both approaches investigated in \cite{OpenQuest}
(see also Remark~\ref{rem:defDiff} below).
Although simple and intuitive, our definition \eqref{eq:def} was not chosen as the natural one in time scales \cite{OpenQuest}.
According to \cite{OpenQuest}, the most natural definition of logarithm  (the first approach of \cite{OpenQuest}) is
\begin{equation*}
L(t) =
\int_1^t \frac{1}{\tau + 2 \mu(\tau)} \Delta\tau \, .
\end{equation*}
Differently from both approaches followed in \cite{OpenQuest},
we show here that \eqref{eq:def} gives a proper logarithm
function on time scales.


\section{Properties}

We begin by proving an analogous relation in time scales
to the equality
\begin{equation*}
\frac{d}{dt} \ln\left(p(t)\right) = \frac{p'(t)}{p(t)} \, .
\end{equation*}

\begin{proposition}
\label{prop}
Assume $p : \mathbb{T} \rightarrow \mathbb{R}^+$ is strictly increasing and $\tilde{\mathbb{T}} := p(\mathbb{T})$ is a time scale. If $p^{\Delta}(t)$ exist for $t \in \mathbb{T}^{\kappa}$, then
\begin{equation}
\label{eq:my:compLog}
\frac{\Delta}{\Delta t}
L_{\tilde{\mathbb{T}}}\left(p(t)\right)
= \frac{p^{\Delta}(t)}{p(t)} \, .
\end{equation}
\end{proposition}

\begin{proof}
For each $t \in \mathbb{T}^{\kappa}$ the chain
rule \cite[Theorem~1.93]{book:ts} asserts that
\begin{equation*}
\left(L_{\tilde{\mathbb{T}}} \circ p\right)^{\Delta}(t)
= \left[\left(L_{\tilde{\mathbb{T}}}^{\tilde{\Delta}}
\circ p\right) p^{\Delta}\right](t)
= L_{\tilde{\mathbb{T}}}^{\tilde{\Delta}}(p(t))
p^{\Delta}(t) \, .
\end{equation*}
From \eqref{eq:trv:propt} we know that
\begin{equation*}
L_{\tilde{\mathbb{T}}}^{\tilde{\Delta}}(p(t))
= \frac{1}{p(t)} \, ,
\end{equation*}
and the desired result follows.
\end{proof}

\begin{remark}
\label{rem:defDiff}
Proposition~\ref{prop} illustrates well the difference between
our definition \eqref{eq:def} and the second approach of \cite{OpenQuest} that defines the logarithm in such a way that
the delta derivative of the logarithm
coincides with the right-hand side
of \eqref{eq:my:compLog} for any rd-continuous and regressive
$p(\cdot)$. Differently from \cite{OpenQuest},
condition \eqref{eq:my:compLog} holds for
a strictly increasing $p(\cdot)$.
This makes a difference with the classical calculus.
While the chain rule of classical calculus has no monotonicity
assumptions, on time scales we need to require $p(\cdot)$ to be
strictly increasing. A chain rule on time scales that does not
require strictly monotone change of variables is possible,
but involves a different notion of time scale \cite{HS:ts}. We confine ourselves here to the well established delta calculus on time scales and to the question posed in \cite{OpenQuest}. The results of the paper are, however, with the necessary changes, easily formulated for the nabla \cite[Chapter~3]{book:ts1}
or alpha \cite{HS:ts} calculus on time scales.
\end{remark}

\begin{example}
Let $\mathbb{T}=\mathbb{N}$ and $p(t)=t^2$.
It follows from Proposition~\ref{prop} that
$$
\left(L_{p(\mathbb{N})}
\circ p\right)^{\Delta}(t)=\frac{t+\sigma(t)}{t^2}
=\frac{2t+1}{t^2}=L_{\mathbb{N}}^{\Delta}(t)\left(2
+L_{\mathbb{N}}^{\Delta}(t)\right) \, .
$$
\end{example}

We now state the main property of our function
$L_{\mathbb{T}}(\cdot)$. It is Theorem~\ref{m:r}
that justify the name \emph{logarithm} to the function
defined by \eqref{eq:def}.

\begin{theorem}
\label{m:r}
Let $t = a b \in \mathbb{T}$, $a \in \mathbb{T}$ with $a > 0$,
and $b \in \mathbb{R}^{+}$. Then,
\begin{equation}
\label{eq:mp}
L_{\mathbb{T}}(a b) =
L_{\mathbb{T}}(a) + L_{\mathbb{T}/a}(b) \, .
\end{equation}
\end{theorem}

\begin{remark}
The $b > 0$ in Theorem~\ref{m:r} does not belong
necessarily to the time scale $\mathbb{T}$. However,
$b$ always belong to $\mathbb{T}/a$ because $a b \in \mathbb{T}$
($s \in \mathbb{T}/a$ if and only if
there exists a $t \in \mathbb{T}$ such that $t = a s$).
\end{remark}

\begin{remark}
When $\mathbb{T} = \mathbb{R}$
one has $\mathbb{T}/a = \mathbb{R}$
for any $a > 0$. Then,
for any $a$ and $b$ positive,
the equality \eqref{eq:mp}
reduces to the classical relation
$\ln(a b) = \ln(a) + \ln(b)$.
\end{remark}

\begin{proof}
By definition,
\begin{equation*}
L_{\mathbb{T}}(a b) =
\int_{1}^{a b} \frac{1}{\tau} \Delta\tau
\end{equation*}
and from basic properties of the delta integral one has
\begin{equation*}
L_{\mathbb{T}}(a b) =
\int_{1}^{a} \frac{1}{\tau} \Delta\tau
+ \int_{a}^{a b} \frac{1}{\tau} \Delta\tau = L_{\mathbb{T}}(a)
+ \int_{a}^{a b} \frac{a}{\tau} \frac{1}{a} \Delta\tau \, .
\end{equation*}
Using the substitution rule for delta integrals
\cite[Theorem~1.98]{book:ts} with $f(\tau) = a/\tau$
and $\nu(\tau) = \tau/a$ (so that $\nu^{\Delta}(\tau) = 1/a$,
$\tilde{\mathbb{T}} = \mathbb{T}/a$, and $\nu^{-1}(s) = a s$),
we arrive to the intended conclusion:
\begin{equation*}
\begin{split}
L_{\mathbb{T}}(a b) &= L_{\mathbb{T}}(a)
+ \int_{\nu(a)}^{\nu(a b)} \left(f \circ \nu^{-1}\right)(s)
\tilde{\Delta}s \\
&= L_{\mathbb{T}}(a)
+ \int_{1}^{b} \frac{a}{a s} \tilde{\Delta}s \\
&= L_{\mathbb{T}}(a) + L_{\tilde{\mathbb{T}}}(b) \, .
\end{split}
\end{equation*}
\end{proof}

\begin{corollary}
\label{cor:dm}
Let $a>0$, $n\in\mathbb{N}$, and $a^k\in\mathbb{T}$ for $k=0,1,\ldots, n$. Then,
\begin{equation}
\label{eq1:dm}
L_{\mathbb{T}}\left(a^n\right)
= \sum_{k=0}^{n-1}L_{\mathbb{T}/a^k}(a) \, .
\end{equation}
\end{corollary}

\begin{proof}
For $n=1$ we have
$L_{\mathbb{T}}\left(a\right)=L_{\mathbb{T}/1}(a)$.
Let $n\in\mathbb{N}$. Assume that \eqref{eq1:dm} is true for any time scale $\mathbb{T}$ and for $n$. Then,
from Theorem~\ref{m:r},
$L_{\mathbb{T}}\left(a^{n+1}\right)
=L_{\mathbb{T}}\left(a\cdot a^n\right)
=L_{\mathbb{T}}\left(a\right)+L_{\mathbb{T}/a}\left(a^n\right)$. Notice that for $k=0,1,\ldots, n-1$ now $a^k\in \mathbb{T}/a$.
Hence, from the inductive assumption for the time scale $\mathbb{T}/a$, we can state that $L_{\mathbb{T}}\left(a^{n+1}\right)
=L_{\mathbb{T}}\left(a\right)
+\sum_{k=0}^{n-1}L_{\left(\mathbb{T}/a\right)/a^k}(a)
= \sum_{k=0}^{n}L_{\mathbb{T}/a^k}(a)$.
\end{proof}

\begin{remark}
Let $a>0$, $n\in \mathbb{N}$.
If $\mathbb{T}=\mathbb{R}$, then $a^n\in\mathbb{R}$
and $\mathbb{R}/a^k=\mathbb{R}$.
Then, $\ln(a^n)=L_{\mathbb{R}}(a^n)
=\sum_{k=0}^{n-1}L_{\mathbb{R}}(a)=n\ln a$.
\end{remark}

\begin{corollary}
\label{cor}
Let $0 < x \in \mathbb{T}$, and $0 < x/y \in \mathbb{T}$. Then,
\begin{equation}
\label{eq:log:frac}
L_{\mathbb{T}}\left(\frac{x}{y}\right)
= L_{\mathbb{T}}(x) - L_{(y\mathbb{T})/x}(y) \, .
\end{equation}
\end{corollary}

\begin{remark}
The positive real $y$ in Corollary~\ref{cor} does not necessarily belong to $\mathbb{T}$. However,
$y \in \mathbb{T}/z$ with $z = x/y$.
\end{remark}

\begin{remark}
If $y = x$, then \eqref{eq:log:frac} gives
$L_{\mathbb{T}}(1) = 0$.
\end{remark}

\begin{remark}
In the particular case $x = 1$ Corollary~\ref{cor} gives
\begin{equation*}
L_{\mathbb{T}}\left(\frac{1}{y}\right)
= - L_{y\mathbb{T}}(y) \, .
\end{equation*}
\end{remark}

\begin{proof}
Let $z = x/y$. Then, $x = y z$ and by Theorem~\ref{m:r}
\begin{equation*}
L_{\mathbb{T}}(x) = L_{\mathbb{T}}(y z)
= L_{\mathbb{T}}(z) + L_{\mathbb{T}/z}(y) \, ,
\end{equation*}
that is, $L_{\mathbb{T}}(z) = L_{\mathbb{T}}(x)
- L_{\mathbb{T}/z}(y)$.
\end{proof}

\begin{example}\mbox{}
\begin{enumerate}

\item[i)]
Let $\mathbb{T}=q^{\mathbb{N}_0}=\{q^n, n\in \mathbb{N}_0\}$, $q>1$. Then,
\begin{equation*}
L_{\mathbb{T}}(t)=(q-1)\frac{\log t}{\log q} \, , \quad
L_{\mathbb{T}}\left(\mathbb{T}\right)=(q-1)\mathbb{N}_0 \, .
\end{equation*}

\item[ii)]
Let $\mathbb{T}=\overline{q^{\mathbb{Z}}}$, $q>1$. Then,
$L_{\mathbb{T}}(t)=\frac{(q-1)\log t}{\log q}$. Moreover, $L_{\overline{q^{\mathbb{Z}}}}\left(q^{\mathbb{Z}}\right)=(q-1)\mathbb{Z}$. For the particular case $q=2$, $L_{\mathbb{T}}(\cdot)$
gives a map from $2^{\mathbb{Z}}$ to $\mathbb{Z}$.

\item[iii)] For $\mathbb{T}=\mathbb{N}$
the logarithm is an harmonic function
(see \cite[Example~1.45]{book:ts}):
\begin{gather*}
L_{\mathbb{N}}(1) = H_{0} = 0 \, , \quad
L_{\mathbb{N}}(n)=\sum_{k=1}^{n-1}\frac{1}{k}=H_{n-1} \, , \ \ n > 1 \, , \\
L_{\mathbb{N}}\left(\mathbb{N}\right)
=\{H_n: n\in \mathbb{N}_0\} \, .
\end{gather*}
\end{enumerate}
\end{example}


\section{Concept of convexity/concavity}

The fact that functions $L_{\mathbb{T}}(\cdot)$ are continuous and increasing is a direct consequence of the definition.
In the standard case $\mathbb{T}=\mathbb{R}$, another important
property of the natural logarithm is that it is a concave
function. To prove that we have now a similar property,
one needs first to define what we mean
by convexity and concavity of a function
on an interval of a generic time scale $\mathbb{T}$.
As before, let $\mathbb{T}$ be a time scale consisting
of at least two points.

\begin{definition}
\label{def:concave}
Let $I$ be an interval in $\mathbb{R}$ such that the set
$I_{\mathbb{T}}:=I\cap \mathbb{T}$ is a nonempty subset of $\mathbb{T}$.
A function $f$ defined and continuous on
$I_{\mathbb{T}}$ is called
\emph{convex on $I_{\mathbb{T}}$}
if for any $t_1, t_2\in I_{\mathbb{T}}$
\begin{equation}
\label{eq:dor:convex}
(t_2-t)f(t_1)+(t_1-t_2)f(t)+(t-t_1)f(t_2)\geq 0 \, ,
\quad t\in I_{\mathbb{T}} \, .
\end{equation}
Similarly, $f$ is said to be
\emph{concave on $I_{\mathbb{T}}$} if
\eqref{eq:dor:convex} holds with $\geq 0$
substituted by $\leq 0$.
\end{definition}

\begin{remark}
Note that $I_{\mathbb{T}} \subset \mathbb{T}$
is closed or open, finite or infinite.
\end{remark}

\begin{remark}
For $\mathbb{T}=\mathbb{R}$ Definition~\ref{def:concave} agrees with the standard definition of convexity and concavity of a function. For an arbitrary $\mathbb{T}$ it may happen that  $I_{\mathbb{T}}$ consists only of one or two points. Then,
$f$ is convex and concave on such $I_{\mathbb{T}}$.
Let $t_1, t_2\in I_{\mathbb{T}}$ and $t_1<t_2$.  If $t\in I_{\mathbb{T}}$ and $t_1\leq t\leq t_2$, then $t=\alpha t_1+(1-\alpha)t_2$ with $\alpha=\frac{t_2-t}{t_2-t_1}$, $1-\alpha=\frac{t-t_1}{t_2-t_1}$.
Thus, we define convexity using convex combinations
of points from $I_{\mathbb{T}}$. Indeed, then the condition~\eqref{eq:dor:convex}
can be rewritten as $f(t)=f(\alpha t_1+(1-\alpha)t_2)\leq \alpha f(t_1)+(1-\alpha)f(t_2)$.
\end{remark}

\begin{theorem}
\label{thm:conv}
Let function $f$ be defined on $I_{\mathbb{T}}:=I\cap \mathbb{T}$ and $\Delta$--differentiable on
$I_{\mathbb{T}}^{\kappa}$. If $f^{\Delta}$ is nondecreasing (nonincreasing) on $ I_{\mathbb{T}}^{\kappa}$, then $f$ is convex (concave) on $I_{\mathbb{T}}$.
\end{theorem}

\begin{proof}
If $I_{\mathbb{T}}$ consists of less than three points, then $f$ is simultaneously convex and concave.
Let us assume that $t\in I_{\mathbb{T}}$ and $t_1<t<t_2$. Then,
the condition \eqref{eq:dor:convex} can be rewritten
in the following equivalent form:
\begin{equation}
\label{eq:dor:con2}
\frac{f(t)-f(t_1)}{t-t_1}\leq \frac{f(t_2)-f(t)}{t_2-t} \,.
\end{equation}
Let us prove this relation.
From the mean value theorem \cite[Theorem~1.14]{book:ts1} we have the existence of points $\tau_1, \xi_1\in [t_1,t)$ and $\tau_2,\xi_2\in [t,t_2)$
such that
\begin{equation}
f^{\Delta}(\tau_1)\leq \frac{f(t)-f(t_1)}{t-t_1}\leq f^{\Delta}(\xi_1), \ \mbox{and} \ f^{\Delta}(\tau_2)\leq\frac{f(t_2)-f(t)}{t_2-t}\leq f^{\Delta}(\xi_2)\,.
\end{equation}
As $t_1\leq \xi_1<\tau_2$,
inequality \eqref{eq:dor:con2} holds
from the assumption $f^{\Delta}(\xi_1)\leq f^{\Delta}(\tau_2)$:
$\frac{f(t)-f(t_1)}{t-t_1}\leq f^{\Delta}(\xi_1)\leq f^{\Delta}(\tau_2)\leq\frac{f(t_2)-f(t)}{t_2-t}$.
\end{proof}

\begin{proposition}
Let $\mathbb{T}$ be a time scale with at least two positive points, one of them being one:
$1 \in \mathbb{T}$, there exists at least one
$t \in \mathbb{T}$ such that $0 < t \ne 1$.
The natural logarithm $L_{\mathbb{T}}(\cdot)$
is concave on $\mathbb{T}\cap (0,+\infty)$.
\end{proposition}

\begin{proof}
Since $L_{\mathbb{T}}^{\Delta}(t)=\frac{1}{t}$ is a decreasing function on $\mathbb{T}\cap (0,+\infty)$, the result follows as a corollary of Theorem~\ref{thm:conv}.
\end{proof}

The following result is a consequence
of Theorem~\ref{thm:conv} and
\cite[Corollary~1.16]{book:ts1}:

\begin{corollary}
Let function $f$ be defined and continuous on an interval $I_{\mathbb{T}}$ and let $f^{\Delta^2}$
exist finite on $I_{\mathbb{T}}^{\kappa^2}$.
Then, $f$ is convex (concave) on $I_{\mathbb{T}}$ if $f^{\Delta^2}(t)\geq 0$ ($f^{\Delta^2}(t)\leq 0$) for
all $t\in I_{\mathbb{T}}^{\kappa^2}$.
\end{corollary}



\end{document}